\magnification = \magstep1
\input amssym.def
\input amssym.tex
\phantom{ }
\vskip 2 true cm

\centerline {\bf STEIN MANIFOLDS AND MULTIPLICITY-FREE}
\centerline {\bf REPRESENTATIONS OF 
COMPACT LIE GROUPS}
\bigskip
\bigskip
\centerline{\bf Dmitri Akhiezer}
\bigskip
\bigskip
\centerline{\bf 1. Introduction}
\bigskip
\noindent
This paper is an outgrowth of a talk, whose aim was to give a survey
of recent results on
multiplicity-free representations of compact Lie groups
arising in complex analysis.
Some proofs are given and some results are new.
For example, spherical complex spaces are defined with respect 
to any (and not just compact) real form of a complex
reductive algebraic group $G$.
Theorem 4.3 is a new characterization of such spaces in terms
of representation theory. 
A
stronger result (Theorem 4.4) was known for compact real forms $K \subset G$. 
In that case
one has the so called complexification theorem which
allows to consider Stein spaces acted on by $K$ as $K$-invariant domains
in Stein spaces with $G$-action.
For a non-compact real form we have to make an 
additional assumption that there exists an open
embedding with similar properties.

In the theory of compact complex homogeneous manifolds, 
a considerable role is played by manifolds fibered by tori
over flag manifolds. It turns out that a compact complex
homogeneous manifold $X$ of a reductive algebraic group is spherical
if and only if $X$ has such a fibration, see Theorem 3.5.
This new result shows that non-algebraic spherical complex    
spaces have interesting geometric properties.

All complex spaces are assumed to be reduced and irreducible.
The central question of the paper is the existence of
antiholomorphic involutions which have nice behaviour with respect to group
actions.
Namely, let $K$ be a
connected compact Lie group acting by holomorphic transformations
on a complex space $X$. For $X$ spherical and Stein,
we want to construct an antiholomorphic involution $\mu : X \to X$
which is $\theta$-equivariant with respect to a Weyl involution 
$\theta : K \to K$.
The reason is that, for holomorphically separable complex spaces,
it is then easy to check whether the algebra ${\cal O}(X)$
is multiplicity-free as a $K$-module.
This is the case if and only if $\mu $ preserves all $K$-orbits,
see Theorem 6.4.
Moreover, for any complex space
the existence of some antiholomorphic involution
preserving each $K$-orbit is sufficient for ${\cal O}(X)$ to 
be multiplicity-free, see Theorem 5.2.

A $\theta $-equivariant antiholomorphic involution $\mu $  
is constructed for spherical Stein 
manifolds, see [A2] and Theorem 6.9 below. 
The key point is the algebraic part of the theorem, i.e., the 
similar construction for smooth affine spherical varieties.  
Any such variety is a fiber bundle over a spherical affine
homogeneous variety $G/H$ with a spherical $H$-module as a fiber.
The isotropy subgroup $H\subset G$ is reductive
and belongs to the class of so called adapted subgroups, see Sect. 6.
We construct $\mu $ with required properties
for $X =G \times _H V$, where $H$ is a connected
adapted subgroup and $V$ is an arbitrary regular $H$-module.
The latter construction is new, see Theorem 6.7. 
This theorem does not solve the problem 
if $H$ is disconnected. The proof in the general
case can be found in 
[A2].

\bigskip
\bigskip
\centerline {\bf 2. Preliminaries}
\bigskip
\noindent
Let $K$ be a compact Lie group. We will use the standard notation
$\hat K$ for the set of all equivalence classes of 
irreducible representations of $K$. For each $\delta \in \hat K$
we denote by $V_\delta $ an irreducible 
$K$-module corresponding to $\delta $. Then, of course,
$V_\delta $ is a finite-dimensional Hilbert space, on which $K$ acts by a
unitary representation. However,
we do not specify an invariant scalar product on $V_\delta $, which is
unique up to a positive multiple. Let $d(\delta )$ be the dimension of 
$V_\delta $,    
$\xi _\delta $ the character of $V_\delta $,
and $\chi _\delta = d(\delta )\xi _\delta $.

We now recall a fundamental result of Harish-Chandra on the
representations of compact Lie groups in Fr\'echet spaces,
see [H].
Let $\varrho$ be a continuous representation of $K$ in a Fr\'echet
space $F$. Then
there is a family of continuous operators
$$E_\delta : F \to F,\ \ E_\delta (f) = \int _K\overline {\chi _\delta (x)}
\cdot \varrho(x)f \cdot d\mu (x) \ \ (f\in F, \ \delta \in \hat K),$$
where $\mu $ is the normalized Haar measure on $K$. 
The orthogonality relations 
for characters imply that
$$E_\delta ^2 = E_\delta , \ \ E_\delta \cdot E_\epsilon =0 
\ \ (\delta \ne \epsilon),$$
and each $E_\delta $ commutes with all $\varrho (k), \  k\in K$.
It follows that $$F_\delta := E_\delta (F) = \{f\in F \, \vert 
\, E_\delta f = f\}$$
is a closed $K$-invariant subspace in $F$.
\bigskip
\noindent{\sl Definition 2.1.} \ $F_\delta $ is called the isotypic component
of $F$ of type $\delta $. A representation $\varrho $ is called
multiplicity-free if each non-zero isotypic component is irreducible, i.e.,
either $F_\delta =\{0\}$ or $F_\delta $ is isomorphic to $V_\delta $ as a
$K$-module.
\bigskip
\noindent
A vector $f \in F$ is said to be differentiable if the mapping
$K \to F,\ \ k \mapsto \varrho (k)f$, is of class $C^\infty $. The subspace
of differentiable vectors is denoted by $F^\infty $. A theorem of G\aa rding 
says that 
$F ^\infty$ is dense in $F$. This result is valid for a continuous 
representation 
of any Lie group countable at infinity. 

A vector $f \in F$ is said to be
$K$-finite if the linear span of the orbit $\varrho (K) f$ is
a finite-dimensional subspace. The subspace $F^0$ of $K$-finite vectors
is contained in $F^\infty $. In our case $K$ is compact and 
we have 
$F_\delta \subset F^0$ for all $\delta \in \hat K$. Indeed, for 
$f \in F_\delta $
and $k \in K$ one has
$$ \varrho (k)f = \varrho (k)E_\delta f = E_\delta \varrho (k) f =
\int _K \overline {\chi _\delta (x)}\cdot \varrho(xk)f \cdot d\mu(x)
= \int _K \overline {\chi _\delta (xk^{-1})}\cdot \varrho(x)f \cdot d\mu(x)
$$
by the invariance of $\mu $.
Now, if $a_{ij}(x)$ are the matrix elements with respect to some basis
of $V_\delta $ and the vectors $f_{ij} \in F$ 
are defined by
$$f_{ij} =\int _K \overline {a_{ij}(x)}  \cdot \varrho(x)f \cdot d\mu (x) 
$$
then $\varrho (K) f $ is contained in the 
linear span of $f_{ij}$.

\bigskip\noindent{{\bf Theorem 2.2 }\ (Harish-Chandra, see [H]). 
{\it For $f \in F^\infty $
one has
$$ f = \sum _{\delta \in \hat K}\ E_\delta f\ ,$$
where the convergence is absolute with respect to any continuous 
seminorm on $F$. The series is finite (has only finitely many non-zero terms)
if and only if $f \in F^0$.}

\bigskip
\noindent
This is a central result which has many applications. For example, 
Theorem 2.2 shows that $F^0$ is dense in $F$. 
Here is another corollary. Recall that the commuting algebra
of a representation $\varrho $ is the algebra of all continuous
endomorphisms $A:F \to F$ 
satisfying $A\varrho (k) = \varrho (k)A$ for all $k \in K$.

\bigskip \noindent{\bf Theorem 2.3.}\ {\it A 
representation $\varrho $ is multiplicity-free if and only if the commuting
algebra ${\cal A}(\varrho )$
is commutative.               }

\bigskip
\noindent
{\it Proof.}\
An endomorphism from ${\cal A}(\varrho)$ preserves the isotypic 
components.
If $\varrho $ is multiplicity-free then, by Schur's lemma,
two such endomorphisms
commute on each $F_\delta $.
By Theorem 2.2 they commute everywhere on $F$.

Conversely, assume $\varrho $ is not multiplicity-free. Then
at least one $F_\delta $ contains two distinct irreducible submodules
$V_1$ and $V_2$. On $V=V_1 + V_2$, we can find two endomorphisms $A$ and $B$
commmuting with the group action on $V$ with $AB \ne BA$. By Hahn-Banach 
theorem there exists a closed linear subspace $W \subset F$, such 
that $F = V \oplus W$. Extend $A $ and $B$ to $F$ by $A\vert W = B\vert W = 0$
and define
$$\tilde A = \int _K \varrho (k) A \varrho (k)^{-1} \cdot d\mu (k),\ \ 
   \tilde B = \int _K \ \varrho (k) B \varrho (k)^{-1} \cdot d\mu (k).$$
Then $\tilde A, \tilde B \in {\cal A}(\varrho )$ and 
$\tilde A\tilde B \ne \tilde B \tilde A$ by construction. 
\hfill $\square $

\bigskip
\noindent{\it Remark 2.4.}\
Multiplicity-free unitary representations (of arbitrary groups)
are usually defined as those for which the commuting algebra
is commutative. We refer the reader to [Ko], see Definition 1.5.3, 
where the notion  
is generalized to representations in topological vector spaces. 
In our context, i.e.,
for compact Lie groups acting in Fr\'echet spaces, this amounts
to Definition 2.1. 

\bigskip
\noindent
We are mostly interested in representations arising from transformation groups
of complex spaces. So let $X$ be a reduced complex space, ${\cal O }(X)$
the corresponding algebra of holomorphic functions, and 
$K$ a Lie group acting on $X$ by holomorphic transformations. 
Then the action
of $K$ on $X$ induces a linear representation
in the Fr\'echet space $F = {\cal O}(X)$, namely,
$$\varrho (k)f (x) = f(k^{-1}x) \ \ (k\in K, x\in X),$$
and we write ${\cal O}_\delta (X)$ for $F_\delta $.
In this setting we have $F^\infty = F$, see e.g. [A1], Sect. 5.2.
Thus, if $K$ is compact we obtain the following corollary of Theorem 2.2.

\bigskip
\noindent{\bf Theorem 2.5.}\ {\it For any $f \in {\cal O}(X)$ one has
$$f = \sum _{\delta \in \hat K}\ E_\delta f\ ,$$
where the series is compactly convergent.
The series is finite if and only if $f$ is $K$-finite.}

\bigskip
\noindent{\it Example 2.6.}\ A
Reinhardt domain $D \subset {\Bbb C}^n$ 
is a domain invariant under the action
of $K =  {\rm {U}}(1)^n$ by rotations
$$(z_1, \ldots , z_n) \to (e^{i\varphi _1}z_1, \ldots, e^{i\varphi _n}z_n).$$
The series in Theorem 2.5 is the Laurent series of $f \in {\cal O}(D)$.
All non-zero isotypic components of $F ={\cal O}(D)$ are one-dimensional,
and so ${\cal O}(D)$ is multiplicity-free. 

\bigskip
\noindent
{\it Example 2.7.}\
Let $V$ be a complex vector space,
$K \subset {\rm {GL}}(V)$ a connected compact linear group,
and $G = K^{\Bbb C}\subset {\rm GL}(V)$ the complex reductive group
containing $K$ as a maximal compact subgroup. Assume that
\smallskip
$(*\,)$ each irreducible $G$-module
occurs in ${\Bbb C}[V]$ at most once. 
\smallskip
\noindent
Then the expansion of  
$f \in {\cal O}(V)^0$ in a series of homogeneous polynomials has
only finitely many non-zero terms,
for otherwise the series from Theorem 2.5 would be infinite. 
In other words, $f \in {\Bbb C}[V]$.
Moreover,
an isotypic component of ${\cal O}(V)$ is contained in
a subspace of homogeneous polynomials and is in fact
irreducible. Thus $(*\,)$ holds if and only if the representation of $K$ in
${\cal O}(V)$ is multiplicity-free.  

We have shown that $(*\, )$ yields the equality
${\cal O}(V)^0 = {\Bbb C}[V]$.
Furthermore, if 
$D \subset V$ is a $K$-invariant domain containing $0\in V$ then
$(*\, )$ implies that ${\cal O}(D)^0 = {\Bbb C}[V]$. To see this, 
introduce a $K$-invariant Hermitian metric in $V$, take a ball in $D$
centered at 0, restrict $f \in {\cal O}(D)^0$ to
that ball and apply the above argument.
A complete list of irreducible 
reductive linear groups with property $(*\,)$
was found by V.Kac, see [K]. The result was extended by A.S.Leahy  
to all connected reductive linear groups, see [L].

\bigskip
\noindent
{\it Example 2.8.}\
Let $D$ be a bounded symmetric domain and
let $K$ be the isotropy subgroup of a point $o\in D$. 
Then the representation of $K$ in ${\cal O}(D)$ is multiplicity-free.

One proof of this fact can be obtained from Theorem 4.4. Namely, 
the isotropy representation of $K$ in the holomorphic
tangent space $V = T_o (D)$ gives rise to a multiplicity-free 
polynomial algebra ${\Bbb C}[V]$, see [J]. Of course, if we assume
$D$ irreducible then $V$ as a $K^{\Bbb C}$-module
occurs in the list in [K]. 
It is known that $D$ can be realized as a $K$-invariant domain in $V$.
Therefore the representation of $K$ in ${\cal O}(D)$ 
is multiplicity-free by
Theorem 4.4, ${\rm (c)} \Rightarrow {\rm (b)}$. 
For another proof see Example 5.4.

\bigskip \noindent
In the sequel, we will always assume that the compact Lie group $K$
is connected. 
This is justified by the following theorem.

\bigskip \noindent{\bf Theorem 2.9.}\ {\it Let 
$K$ be a compact Lie transformation group of an
irreducible reduced complex space $X$. Let $K^\circ \subset K$ be the connected   
component of the neutral element. Then ${\cal O}(X)$ is
multiplicity-free with respect to $K$ if and only if ${\cal O}(X)$
is multiplicity-free with respect to $K^\circ $.}

\bigskip
\noindent
{\it Proof.}\
If ${\cal O}(X) $ is multiplicity-free under $K^{\circ }$
then the same
is obviously true for $K$. To prove the converse, assume $V_1, V_2$
are two distinct isomorphic $K^\circ $-modules in ${\cal O}(X)$.
Let $f_1 \in V_1,\ f_2 \in V_2$ 
be their highest weight vectors with respect to some Borel subalgebra 
of the complexified Lie algebra of $K^\circ $. 
Note that for any integer $m$ the functions 
$f_1^k \cdot f_2^{m-k},\ k=0, \ldots, m$, are linearly independent 
highest weight vectors of the same weight. 
Thus there exists $\delta _0 \in \widehat {K^\circ }$,
such that  
$${\rm dim}\, {\cal O}_{{\delta }_0}(X) > N\cdot d(\delta _0),$$ 
where $N>0$ is an arbitrary constant. Take
$N = [K:K^\circ ]$ and let $\{\delta _0, \delta _1, \ldots, \delta _p\}
\subset \widehat {K^\circ }$
be the set of all equivalence classes of irreducible representations obtained
from $\delta _0$ by an inner automorphism of $K$. Let $V$
denote the $K$-submodule of ${\cal O}(X)$ generated by 
${\cal O}_{\delta _0}(X)$.
Then
$V$ is the direct sum of ${\cal O}_{\delta _i}(X)$,
$\ i=0,1,\ldots , p$, and all these subspaces have the same dimension.
Therefore
$${\rm dim}\,V = (p+1)\cdot {\rm dim}
\,{\cal O}_{\delta_0}(X) > (p+1)Nd(\delta _0).$$
Let $[\delta : \delta _i]$ be the number of times
$\delta _i$ occurs in the restriction of $\delta \in \hat K$ to $K^\circ $.
Frobenius reciprocity theorem shows that
$$\sum _{\delta \in \hat K} [\delta : \delta _i]\cdot d(\delta ) =
Nd(\delta _i) = Nd(\delta _0)$$
for all $i$, $i=0,1, \ldots, p$. 
Now, $V$ is a multiplicity-free $K$-module whose irreducible components
are among $V_\delta $ with $[\delta : \delta _i] \ge 1$ for some $i$.
It follows that
$${\rm dim}\, V \le N(p+1)d(\delta _0),$$
which contradicts the above lower bound.
\hfill $\square$

\bigskip
\bigskip
\centerline {\bf  3. Spherical varieties and their holomorphic counterparts}
\bigskip \noindent
We recall the following definition from the theory of algebraic
groups. All varieties and morphisms
are defined over $\Bbb C$.

\bigskip
\noindent
{\sl Definition 3.1.}\
Let $G$ be a connected reductive algebraic group
acting on a normal algebraic variety $Y$ and let $B \subset G$ be a Borel
subgroup. 
Then
$Y$ is called spherical (with respect to $G$) if $B$ has an open orbit on $Y$.

\bigskip
\noindent
We refer the reader to [Lu], Sect. 0.2, for several 
possible answers to
the question: "Pourquoi s'int\'eresser aux vari\'et\'es sph\'eriques?"
Our motivation comes from the following result.
\bigskip
\noindent
{\bf Theorem 3.2} (see [S], [VK]).
{\it Let $L$ be an algebraic line bundle on $Y$ with $G$-linearization.
If $Y$ is spherical then 
the $G$-module of regular global sections $\Gamma (Y,L)$
is multiplicity-free, i.e., each irreducible
$G $-module occurs in $\Gamma (Y,L)$ with multiplicity 0 or 1. 
If $Y$ is affine and ${\Bbb C}[Y]$
is a multiplicity-free $G$-module then $Y$ is a spherical $G$-variety. 
In particular, a normal affine $G$-variety $Y$ is spherical if and only
if the algebra ${\Bbb C}[Y]$ is a multiplicity-free $G$-module.}
\bigskip
\noindent
In complex analysis, it is quite often that a real form $G_0 $
of a complex Lie group $G$
acts on a complex space $X$ by holomorphic transformations, whereas $G$ 
acts on $X$ only locally. In this setting, the Lie homomorphism makes  
the complex Lie algebra $\goth g$ into 
an algebra of holomorphic vector fields on $X$. A typical example 
is a $G_0$-invariant domain in a complex space acted on by $G$.

\bigskip
\noindent
{\sl Definition 3.3.}\ Let 
$G_0$ be a real form of a connected complex reductive algebraic group $G$.
Let $X$ be a normal (irreducible,
reduced) complex space on which $G_0$ acts by holomorphic 
transformations.
Then $X$ is said to be spherical (with respect to $G_0$)
if there exists a point $x\in X$ such that the holomorphic tangent space
$T_x(X)$ is generated by vector fields from $\goth b$, the Lie algebra of $B$.
\bigskip
\noindent
{\it Remark 3.4.}\
If $G$ acts holomorphically on $X$ then this definition is independent of $G_0$  
and amounts to saying that $B$ has an open orbit on $X$.
The following theorem shows that there
exist non-algebraic actions of this type. 
In what follows, we denote by $L^\prime $ the commutator subgroup of a 
group $L$.
\bigskip
\noindent
{\bf Theorem 3.5.}\ {\it Let $X = G/H$
be a compact complex homogeneous space. Then $X$
is spherical in the sense of Definition 3.3 if and only if $X$ is
a locally trivial holomorphic fiber bundle over a flag manifold with
a complex torus as a fiber.} 

\bigskip
\noindent
{\it Proof.}\ By the normalizer theorem due to
J.Tits, A.Borel and R.Remmert (see e.g. [A1], Sect. 3.5., and references
therein),
the normalizer of the connected component $H^\circ \subset H$ is                 
a parabolic subgroup in $G$, which will be
denoted by $P$. 

Assume first that $X$ is spherical and fix a Borel subgroup $B \subset G$
such that, on the Lie algebra level, $\goth g = \goth b + \goth h$. 
Choose a maximal algebraic torus $T$ in $P\cap B$ and
introduce an ordering of the root system
in such a way that $\goth b = \goth b^-$.
Since $H^\circ $ is normalized by $T$, the Lie algebra $\goth h$  
is spanned by $\goth h \cap \goth t$ and the root vectors $e_\alpha $
which are contained in 
$\goth h$. The decomposition $\goth g = \goth b^- + \goth h$
shows that all positive root vectors are in $\goth h$, hence
$[\goth b^+ , \goth b^+] \subset \goth h \subset \goth p$. Now,
if $\alpha $ is a positive root such that 
$e_{-\alpha} \in \goth p$ then 
$$[e_{-\alpha}, e_\alpha] \in [\goth p, \goth h] \subset \goth h$$
and
$$e_{-\alpha } \in {\Bbb C}\cdot[e_{-\alpha}, [e_{-\alpha }, e_\alpha]]\in
                              [\goth p, \goth h] \subset \goth h.$$
Therefore $[\goth p, \goth p] \subset \goth h$ and $P^\prime \subset H $,
showing that the fiber $P/H$ of the fibration $G/H \to G/P$ is a 
complex torus.

To prove the converse, we recall that 
compact complex homogeneous manifolds, fibered by tori
over flag manifolds, are well-understood. In particular, 
the existence of a parabolic subgroup $P \subset G$ with $P^\prime \subset H
\subset P$ follows from the normalizer theorem, see [A1], Sect. 3.5,
Cor. 2, and Sect. 3.6, Prop. 1.
On the other hand,
if $P$ is a parabolic subgroup then
$G/P^\prime $ is a spherical algebraic variety. Thus, for any closed
complex Lie subgroup $H$ containing $P^\prime $, the corresponding homogeneous
space $G/H$ is spherical in the sense of Definition 3.3.
\hfill $\square $
\bigskip
\noindent
{\bf Corollary 3.6.}\ {\it Let $X = G/H$
be a compact complex homogeneous space. Then $X$  
spherical in the sense of Definition 3.3 if and only if
there exists
a parabolic subgroup $P \subset G$, such that $P^\prime \subset H \subset P$
and $P/H$ is compact. Moreover, such an $X$ is an algebraic 
homogeneous space (which is then spherical 
in the sense of Definition 3.1) if and only if $H = P$ and $X$ is a flag
manifold.}

\bigskip \noindent
{\it Proof.}\ The first statement follows from the proof of Theorem 3.5.
The second statement is an extract from the theory of
algebraic groups.
Namely, if $X= G/H$ is an algebraic homogeneous space then $H$ 
is an algebraic subgroup of $G$. By Chevalley's theorem $X$
is a quasi-projective variety with an equivariant projective embedding. 
Since $X$ is also compact,
the image of this embedding is closed 
and Borel's fixed point theorem shows that $H$
is parabolic.
\hfill $\square$         
 
\bigskip
\bigskip
\centerline {\bf 4. Spherical Stein spaces }

\bigskip
\noindent
In this section, we carry over Theorem 3.2 to the category
of complex spaces. 
\bigskip
\noindent
{\bf Lemma 4.1} (see [AH], Lemma 2).\ 
{\it
Let $\Omega \subset {\Bbb C}^n $ be a domain containing the origin and let
$$A_i = \sum _{j=1}^n a_{ij}(z)\,{\partial \over {\partial z_j}}$$
be holomorphic vector fields in $\Omega $, such that $a_{ij}(0) = \delta_{ij}$.
If $f \in {\cal O}(\Omega )$, 
$A_if\in {\cal O}(\Omega )f$, and $f(0) = 0$ 
then $f=0$.}
\bigskip
\noindent
{\bf Theorem 4.2 } (cf. [AH], Theorem 1).  
\ {\it We use the notations from Definition 3.3.
Let $X$ be a spherical complex space with respect to $G_0$
and let $L$ be a holomorphic line bundle on $X$ with $G_0$-linearization.
Then each irreducible finite-dimensional $G_0$-module occurs
in $\Gamma (X,L)$ with multiplicity 0 or 1.}

\bigskip
\noindent
{\it Proof.}\ We can find a non-singular point $x\in X$, 
a coordinate neighborhood $\Omega $ of 
$x \in X$ with coordinates $z_i,\ z_i(x)=0,$
and the vector fields $A_i$ in $\Omega $ coming from $\goth b$ and 
satisfying the assumptions of Lemma 4.1.
Let $V_1$ and $V_2$ be two distinct, but isomorphic irreducible 
finite-dimensional $G_0$-submodules
in $\Gamma (X,L)$. 
Note that $V_1$ and $V_2$ are also irreducible and isomorphic 
$\goth g$-modules.
Thus they have the same highest weight $\lambda: \goth b \to 
{\Bbb C}$, and so we obtain
two linearly independent holomorphic global sections $s_i \in
\Gamma (X,L)$, 
such that
               $$As_i = \lambda (A)s_i\ \ \ (A\in \goth b,\ i=1,2).$$
Let $s$ be a non-zero linear combination of $s_1, s_2$ vanishing at $x$.
Shrinking $\Omega $, we can find
a section $s_0 \in \Gamma (\Omega, L)$ without zeros. 
Define $\varphi _A \in {\cal O}(\Omega )$ by $As_0 = \varphi _As_0$
and write
$s\vert \Omega = fs_0$ for some $f\in {\cal O}(\Omega )$.
Then
$$Af = (\lambda (A) - \varphi _A)f \in {\cal O}(\Omega )f$$
by Leibniz rule. Therefore 
$f=0$ by Lemma 4.1. Since $s \ne 0$,
we get a contradiction.
\hfill $\square $

\bigskip
\noindent               
{\bf Theorem 4.3.}\ {\it Let $Y$ be a 
normal Stein space with holomorphic action of $G$ and let
$X$ be a $G_0$-invariant domain in $Y$. The following conditions
are equivalent:
\smallskip \noindent 
{\rm (a)} $X$ is spherical with respect to $G_0$;
\smallskip \noindent
{\rm (b)} any irreducible 
finite-dimensional $G_0$-module occurs in ${\cal O}(X)$                     
with multiplicity 0 or 1;
\smallskip \noindent
{\rm (c)} $Y$ is a spherical affine variety of $G$.
}
\bigskip
\noindent
{\it Proof.}\
${\rm (a) }\Rightarrow {\rm (b)}$
by Theorem 4.2 and  ${\rm (c)} \Rightarrow {\rm (a) } $
is
evident. To prove ${\rm (b)}\Rightarrow {\rm (c)}$, 
observe that the trivial
$G$-module occurs in ${\cal O}(Y)$ only once. For, if 
$f_1, f_2 \in {\cal O}(Y)$ were $G$-invariant and non-proportional then
$f_1\vert X, f_2\vert X$ were
$G_0$-invariant and also non-proportional, contradictory to (b).
Since
two closed 
$G$-orbits on $Y$ are separated by $G$-invariant holomorphic functions,
there is only one such orbit. It follows that $Y$ is
($G$-equivariantly biholomorphic to) an affine algebraic variety
on which $G$ acts algebraically, see [Sn], Cor. 5.6. 
Furthermore, ${\Bbb C}[Y]$ is a multiplicity-free
$G$-module. Indeed, two distinct isomorphic irreducible $G$-submodules
of ${\Bbb C}[Y]$ would induce by restriction to $X$ two distinct
isomorphic irreducible $G_0$-submodules
of ${\cal O}(X)$, which again contradicts (b). By Theorem 3.2
the variety $Y$ is spherical.
\hfill $\square $

\bigskip
\noindent
For $G_0$ compact one has a stronger result.
\bigskip
\noindent
{\bf Theorem 4.4} (see [AH], Theorem 2).
\ {\it
Let $K$ be a compact real form of $G$ and let $X$ be a normal Stein $K$-space.
Then the following conditions are equivalent:
\smallskip \noindent
{\rm (a)} $X$ is spherical with respect to $K$;
\smallskip \noindent
{\rm (b)} ${\cal O}(X)$ is a multiplicity-free $K$-module;
\smallskip \noindent
{\rm (c)} $X$ is a $K$-invariant domain in a spherical affine $G$-variety.}

  \bigskip
\noindent
{\it Proof.}\ Again, in view of Theorem 4.2 we only have to prove ${\rm (b)}
\Rightarrow {\rm (c)}$. Now,
for compact transformation groups one has the so called
complexification theorem, see [He]. Namely, there exists 
another reduced
Stein space $X^{\Bbb C}$ acted on by $G$
along with a $K$-equivariant open embedding $i: X \hookrightarrow
X^{\Bbb C}$. Moreover, 
$G\cdot i(X) = X^{\Bbb C}$, so $X^{\Bbb C}$ is normal if $X$ is normal.
To finish the proof, it remains to take $Y = X^{\Bbb C}$ in Theorem 4.3.
\hfill $\square $

\bigskip
\bigskip
\centerline {\bf 5. Antiholomorphic involutions}

\bigskip \noindent
For complex manifolds, J.Faraut and E.G.F.Thomas gave an interesting
and simple geometric condition which implies that the function algebra
is a multiplicity-free module, see [FT]. In their setting, 
the group of holomorphic transformations acting on the manifold $X$  
is not necessarily compact and the condition guarantees that any
invariant Hilbert subspace in ${\cal O}(X)$
is multiplicity-free (see Remark 2.4).
We state their result for compact Lie groups acting on complex 
spaces.
 
\bigskip
\noindent
{\sl Definition 5.1.}\
An antiholomorphic self-map $\mu $ of a complex space $X$ is called an
antiholomorphic involution if $\mu ^2 = {\rm {id}}$.

\bigskip
\noindent
{\bf Theorem 5.2 } (see [FT], Theorem 3). 
\ {\it Let $X$ be an irreducible reduced complex space and $K$ a compact
Lie group acting on $X$ by holomorphic transformations. Assume 
there exists an antiholomorphic involution $\mu : X \to X$ such that 
$\mu (x) \in K\cdot x$ for every $x \in X$. Then ${\cal O}(X)$
is a multiplicity-free $K$-module.}

\bigskip
\noindent
{\it Remark 5.3.}\ 
The proof in [FT] goes without changes for irreducible reduced 
complex spaces. A simplified proof is given in [AP], see Proposition 3.3.

\bigskip
\noindent
{\it Example 5.4.}\
Let $D$ be a bounded symmetric domain in a complex vector space $V$.
Further, let $\goth g$ be the Lie algebra 
of the automorphism group ${\rm {Aut}}(D)$ 
and let
$\goth k$ be the subalgebra of $\goth g$ associated to a maximal compact
subgroup $K \subset {\rm {Aut} }(D)$.
We assume that $0 \in D$ and $K$ is the isotropy subgroup of $0$.
In the Cartan decomposition $\goth g = \goth k + \goth p$, we  
can identify $\goth p$ with $V$ (as a real vector space) and consider
a
maximal abelian subspace $\goth a \subset \goth p$ as a real vector 
subspace of $V$.                                              

For every bounded symmetric domain 
there exists an antilinear map $\mu : V \to V$, such that $\mu ^2 = 
{\rm {id}}$,
$\mu (D) = D$ and $\mu (z) = z $ for all $z$ in some Cartan 
subspace $\goth a$. For irreducible bounded symmetric domains
this fact is checked 
using the classification (see [FT]), and the general case follows easily.

Now, any $x \in D$ can be written as $x =k\cdot z$ with $k\in K,\, 
z \in \goth a \cap D$,
and so we obtain $\mu (x) = \mu (k\cdot z) = \mu k \mu (z) \in K\cdot x$.
Thus, Theorem 5.2 implies that ${\cal O}(D)$ is a multiplicity-free $K$-module.

\bigskip
\bigskip
\centerline {\bf 6. Spherical Stein manifolds and the Weyl involution }
\bigskip
\noindent
In the final section, we discuss the converse to Theorem 5.2.  
For this, additional assumptions on $X$ are necessary. Namely, 
$X$ must have sufficiently many holomorphic functions. 
Indeed, if e.g. ${\cal O}(X) = {\Bbb C}$ and $K$ is trivial then ${\cal O}(X)$
is multiplicity-free, but a self-map of $X$ preserving $K$-orbits
is the identity map which is holomorphic and not antiholomorphic.
In the end, we will assume that $X$ is a Stein manifold, but
our first result holds in a more general setting.
To state this result, we need several definitions.

\bigskip
\noindent
{\sl Definition 6.1.}\ Given a group $\Gamma $ acting on two sets
$A$ and $B$ and an automorphism $\vartheta: \Gamma \to \Gamma$, one says
that a map
$\mu: A \to B$ is $\vartheta $-equivariant
if $\mu (\gamma \cdot a) = \vartheta (\gamma) \cdot \mu (a)$ for all 
$a\in A,\, \gamma \in \Gamma$.

\bigskip
\noindent
{\sl Definition 6.2.}\ An involutive automorphism $\theta $
of a connected compact Lie group $K$ is called a Weyl involution 
if
$\theta (t) = t^{-1}$ for all $t$ in a maximal torus $T \subset K$.

\bigskip
\noindent
{\it Remark 6.3.}\ It 
is known that such an involution exists and that two Weyl
involutions are conjugate by an inner automorphism of $K$, see e.g. 
[W], Sect. 12.6. 
From the point of view of the representation theory,
the main property of $\theta $ is the following one.
For a representation $
\varrho:K \to {\rm {GL}}(V)$, the representation $k \mapsto 
\varrho(\theta (k))$
in the same vector space $V$ is dual to the given one.
The notion of a Weyl involution $\theta $
can also be defined for a connected reductive algebraic group $G$ over $\Bbb C$.
The definition and the properties of $\theta : G \to G$ are similar.

\bigskip
\noindent
{\bf Theorem 6.4 } (see [AP], Theorem 4.1).
\ {\it Let $X$ be a holomorphically separable irreducible reduced complex
space, $K$ a connected compact Lie group acting on $X$ by holomorphic
transformations, $\theta : K \to K$ a Weyl involution, and
$\mu : X \to X$ a $\theta$-equivariant antiholomorphic
involution of $X$. Then ${\cal O}(X)$ is a multiplicity-free 
$K$-module if and only if $\mu (x) \in K\cdot x$ for every $x \in X$.
}
\bigskip
\noindent
{\it Proof.}\ Without going into details, we explain
the idea, showing the role of the Weyl involution.
Let $V \subset {\cal O}(X)$ be an irreducible $K$-submodule.
Introduce a $K$-invariant Hermitian inner product and choose
a unitary basis $\{f_j\}$ in $V$. It is easily seen that
the function 
$$\Phi=\sum _j\,  f_j \bar{f_j}$$
is $K$-invariant and does not depend on the choice of basis.
Assuming ${\cal O}(X)$ multiplicity-free, we 
get a family of real-analytic functions $\{\Phi_\delta\}$,
one for each isotypic component ${\cal O}_\delta (X)$. Using
the fact that $X$ is holomorphically separable, one can prove
that these functions separate $K$-orbits.

On the other hand, let $\mu f(x) = f(\mu (x))$ for any function $f$ on $X$.
Since $\mu $ is $\theta $-equivariant, $\mu {\cal O}_\delta (X)$
is an irreducible $K$-module with representation twisted by $\theta $.
But $\theta $ is the Weyl involution, so this $K$-module
is dual to ${\cal O}_\delta (X)$. Since the $K$-module 
$\overline {{\cal O}_\delta (X)}$ is also dual to ${\cal O}_\delta (X)$
and ${\cal O}(X)$ is multiplicity-free, we have the equality 
$$\mu {\cal O}_\delta (X) = \overline {{\cal O}_\delta (X)}.$$
Moreover, one can show that the composition of $\mu $ 
with complex conjugation preserves a $K$-invariant Hermitian inner product
in ${\cal O}_\delta (X)$. From this it follows that $\mu \Phi_\delta 
= \Phi_\delta $.
Since the family $\{\Phi_\delta \}$ separates $K$-orbits, $\mu $ must
preserve each of them, i.e., $\mu (x) \in K\cdot x$ for all $x \in X$.

The converse is already known from Theorem 5.2.
\hfill $\square $

\bigskip
\noindent
{\bf Lemma 6.5 } (see [A2], Theorem 3.1).
\ {\it Let $L$ be a connected compact Lie group, $\theta : L \to L$
a Weyl involution, and $\varrho : L \to {\rm {GL}}(V)$
a
complex representation of $L$. Then there exists an antilinear map 
$\nu:V \to V$, such that $\nu ^2 = {\rm {id}}$ and
                           $$\nu (\varrho(l)v) = \varrho(\theta (l))\,\nu (v)$$
for all $l \in L,\, v\in V$.
If $\varrho $ is irreducible and $\nu ^\prime : V \to V$ is
another map with the same properties then $\nu ^\prime = c\nu$ for some $c$
with $\vert c \vert = 1$.}
\bigskip
\noindent
{\sl Definition 6.6 } (see [AV]). 
\ Let $G$ be a connected complex reductive group
and $H \subset G$ a reductive algebraic subgroup. Then $H$ is called 
adapted if there is a Weyl involution $\theta : G \to G$,
such that $\theta (H) = H$ and $\theta $ induces a Weyl involution
of the connected component $H^\circ $.
\bigskip
\noindent
{\bf Theorem 6.7.}\ {\it Let $H \subset G$ be connected and adapted
and let $X =G\times_H V$ be a homogeneous fiber bundle
with $H$ acting on $V$ by a regular linear representation.
For any maximal compact subgroup $K \subset G$ and any Weyl involution
$\theta : K \to K$ there exists a $\theta $-equivariant
antiholomorphic involution $\mu : X \to X$.} 
\bigskip
\noindent
{\it Proof.}\ It suffices to prove the theorem for some $K$ and $\theta $.
Let $L$ be a maximal compact subgroup of $H$. Note that
$L$ is connected. Choose a maximal compact subgroup $K\subset G$
containing $L$. According to [AV], Prop. 5.14, 
we can find a Weyl involution $\theta $ of $G$ so that $K$ and $L$
are $\theta $-invariant and the restriction of $\theta $ 
to $K$ and to $L$ is a Weyl involution of these groups. 
Moreover, if $\tau : G \to G$ is the Cartan involution
with fixed point subgroup $K$ then $\tau \theta = \theta \tau $.
The product $\sigma = \tau \theta $ is an antiholomoprphic involution of $G$.
We remark that the fixed point subgroup of $\sigma $ is a split real 
form of $G$.

We can identify $H$ and $L$ with their images in ${\rm {GL}}(V)$
and assume that $\varrho $ in
Lemma 6.5 is the identity representation.
By that lemma there exists an antilinear involution $\nu : V \to V$,
such that 
$$\nu (lv) = \theta (l)\nu(v)$$
for all $l \in L,\, v\in V$.
We claim that in fact 
$$\nu (hv) = \sigma (h)\nu(v)$$
for all $h \in H$.
Indeed, if $v \in V$ is fixed then the second equality holds
for all $h$ in a maximal totally real submanifold $L \subset H$
(where it is just the first one) and therefore everywhere on $H$.

Consider the antiholomorphic involution $\tilde \mu $ of $G \times V$,
defined by
$$\tilde \mu (g,v) = (\sigma (g), \nu (v)).$$
For $h \in H$ let $t_h(g,v) =(gh^{-1}, hv).$
The transformations $t_h$
define an action of $H$ on $G\times V$, and $X$ is the geometric quotient 
of that action.
Since $\nu (hv) = \sigma (h)\nu(v)$, it follows that
$\tilde \mu \cdot t_h = t_{\sigma (h)} \cdot \tilde \mu $.
Thus $\tilde \mu $ gives rise to a self-map $\mu : X \to X$.
Clearly, $\mu $ is an antiholomorphic involution.
Since $\theta $ and $\sigma $ coincide on $K$, the map
$\mu $ is $\theta $-equivariant with respect to the $K$-action on $X$.
\hfill $\square $
\bigskip
\noindent
{\bf Lemma 6.8 } (see [KVS], Cor. 2.2).
\ {\it Let $G$ be a connected complex reductive group and let $X$
be a smooth affine spherical variety of $G$. Then 
$X = G\times _H V$,
where $H\subset G$ is a reductive subgroup with $G/H$ spherical and $V$
is a spherical $H$-module.}

\bigskip
\noindent
{\bf Theorem 6.9} (see [A2], Theorem 1.2).
\ {\it
Let $X$ be a Stein manifold acted on by a connected compact Lie group $K$
of holomorphic transformations. 
Let $\theta : K \to K$ be a Weyl involution.
If $X$ is spherical with respect to $K$ then there exists
a $\theta $-equivariant antiholomorphic involution $\mu : X \to X$.
Any such involution preserves $K$-orbits.
}
\bigskip
\noindent
{\it Proof.}\ We write $G$ for the complexified group $K^{\Bbb C}$.
As in the proof of Theorem 4.4, we can reduce the statement
to the algebraic case. Namely, $X$ can be embedded as a $K$-invariant domain
into an affine spherical $G$-variety $X^{\Bbb C}$. Moreover,
the $G$-saturation of $X$ in $X^{\Bbb C}$ is the whole $X^{\Bbb C}$,
and so $X^{\Bbb C}$ is non-singular. 
Assume that the theorem is proved for $X^{\Bbb C}$.
Then we have an antiholomorphic involution 
$\mu : X^{\Bbb C} \to X^{\Bbb C}$, which is $\theta$-equivariant
with respect to the $K$-action. By Theorem 6.4 $\mu $ preserves $K$-orbits 
on $X^{\Bbb C}$. Thus the subset $X \subset X^{\Bbb C}$ is $K$-stable. 
The restriction of $\mu $ to $X$ is the involution we are looking for,
and the second assertion is now obvious.

Let now $X$ be a smooth affine spherical
variety. Then Lemma 6.8 displays $X$ as a vector bundle, $X = G\times _H V$.
Moreover, since $G/H$ is a spherical variety,
the subgroup $H$ is adapted, see [AV], Prop. 5.10.
If $H$ is connected then we can apply Theorem 6.7 and finish the proof.
For the general case see [A2].   

\bigskip
\noindent
{\bf Corollary 6.10.}\ {\it Let $X$, $K$ and $\theta $ be as in Theorem 6.9.
The following properties of the action $K\times X \to X$ are equivalent:
\smallskip \noindent
{\rm (a)} ${\cal O}(X)$ is a multiplicity-free $K$-module;
\smallskip          \noindent
{\rm (b)} $X$ is a spherical $K$-manifold;
\smallskip           \noindent
{\rm (c)} there exists an antiholomorphic involution $\mu : X \to X$
preserving $K$-orbits;
\smallskip                \noindent
{\rm (d)} there exists a $\theta $-equivariant antiholomorphic
involution $\mu : X \to X$ preserving $K$-orbits.}
\bigskip
\noindent
{\it Proof.}\ This results from Theorems 4.4, 5.2, and 6.9.
\hfill $\square $
\bigskip
\noindent
{\it Remark 6.11.}\ The assumption 
that $X$ is non-singular is only important in the proof of ${\rm (b)}
\Rightarrow
{\rm (d)}$. The author does not know whether this is true
for spaces with (normal) singularities.

\bigskip
\bigskip
\centerline {\bf References}
\medskip
\noindent
[A1]
D. Akhiezer,
{\sl Lie group actions in complex analysis,}
Vieweg,
1995

\noindent [A2]
D. Akhiezer,
{\sl Spherical Stein manifolds and the Weyl involution,}
Ann. Inst. Fourier, Grenoble,
59
(2009),
1029 - 1041

\noindent [AH] D. Akhiezer, P. Heinzner,
{\sl Spherical Stein spaces,}
Manuscripta Math.
114
(2004),
327 - 334

\noindent [AP]
D. Akhiezer, A. P\"uttmann,
{\sl Antiholomorphic involutions of spherical complex spaces,}
Proc. Amer. Math. Soc.
136
(2008),
1649 - 1657

\noindent [AV]
D. Akhiezer, E. B. Vinberg,
{\sl Weakly symmetric spaces and spherical varieties,}
Transform. Groups
4
(1999),
3 - 24

\noindent [FT]
J. Faraut, E. G. F. Thomas,
{\sl Invariant Hilbert spaces of holomorphic functions,}
J. Lie Theory
9
(1999),
383 - 402

\noindent [H]
Harish-Chandra,
{\sl Discrete series for semisimple Lie groups II,}
Acta. Math.
116
(1966),
1 - 111

\noindent [He]
P. Heinzner,
{\sl Geometric invariant theory on Stein spaces,}
Math. Ann.
289
(1991),
631 - 662

\noindent [J]
K. Johnson,
{\sl On a ring of invariant polynomials on a Hermitian symmetric space,}
J. Algebra
67
(1980),
72 - 81

\noindent [K]
V. Kac,
{\sl Some remarks on nilpotent orbits,}
J. Algebra
64
(1980), 
190 - 213

\noindent [Ko]
T. Kobayashi,
{\sl Multiplicity-free representations and visible actions on 
complex manifolds,}
Publ. RIMS, Kyoto,
41
(2005),
497 - 549

\noindent [KVS]
F. Knop, B. Van Steirteghem,
{\sl Classification of smooth affine spherical varieties,}
Transform. Groups
11
(2006),
495 - 516

\noindent [L]
A. S. Leahy,
{\sl A classification of multiplicity free representations,}
J. Lie Theory
8
(1998),
367 - 391

\noindent [Lu]
D. Luna,
{\sl Vari\'et\'es sph\'eriques de type A,}
Publ. Math. de l'IH\'ES
94
(2001),
161 - 226

\noindent [S]
F. J. Servedio,
{\sl Prehomogeneous vector spaces and varieties,}
Trans. Amer. Math. Soc.
176
(1973),
421 - 444

\noindent [Sn]
D. M. Snow,
{\sl Reductive group actions on Stein spaces,}
Math. Ann.
259
(1982),
79 - 97

\noindent [VK]
E. B. Vinberg, B. N. Kimel'feld,
{\sl Homogeneous domains in flag manifolds and spherical subgroups of
semisimple Lie groups,}
Funct. Anal. Appl. (Funktsional. Anal. i ego Prilozhen.)
12
(1978),
168 - 174

\noindent [W]                                
J. A. Wolf,                              
{Harmonic analysis on commutative spaces,}    
Amer. Math. Soc.,                                 
Providence, RI,             
2007

\bigskip
\noindent
Institute for Information Transmission Problems,

\noindent
B. Karetny per. 19, Moscow,

\noindent
127994, Russia

\noindent
{\it E-mail address:}\  akhiezer@iitp.ru
\end